\documentclass[11pt,reqno]{amsart}

\usepackage{latexsym}
\usepackage{amsfonts}
\usepackage{amsmath}
\usepackage{url}
\usepackage{graphicx}
\usepackage{color}
\usepackage{bbm}
\definecolor{myblue}{rgb}{0.09,0.32,0.44} 
\usepackage[pdfborder={0 0 0},pdfborderstyle={},colorlinks=true,linkcolor=myblue,citecolor=myblue,urlcolor=blue]{hyperref}

\usepackage{epsfig,psfrag,verbatim,amsfonts}

\usepackage{amssymb,amsmath}

\def\bl{\begin{lemma}}
\def\el{\end{lemma}}
\def\bth{\begin{theorem}}
\def\eth{\end{theorem}}
\def\bc{\begin{corollary}}
\def\ec{\end{corollary}}
\def\bcj{\begin{conjecture}}
\def\ecj{\end{conjecture}}
\def\bpr{\begin{proposition}}
\def\epr{\end{proposition}}
\def\bde{\begin{definition}}
\def\ede{\end{definition}}
\def\E{\mathbb{E}}
\def\H{\mathbb{H}}

\newcommand{\dist}{\mbox{\rm dist}}
\newcommand{\be}{\begin{eqnarray}}
\newcommand{\ee}{\end{eqnarray}}

\newcommand{\Z}{{\mathbb Z}}

\newcommand{\T}{{\mathcal T}}

\newcommand{\A}{{\mathcal A}}

\renewcommand{\and}{\hbox{ {\rm and} }}

\newcommand{\C}{{\mathcal{C}}}

\newtheorem{theorem}{Theorem}[section]
\newtheorem{definition}{Definition}[section]
\newtheorem{lemma}[theorem]{Lemma}

\newtheorem{corollary}[theorem]{Corollary}
\newtheorem{proposition}[theorem]{Proposition}
\newtheorem{conjecture}[theorem]{Conjecture}

\newtheorem*{theorem*}{Theorem}

\theoremstyle{definition}
\numberwithin{equation}{section}

\begin{document}
\title{Self avoiding walk on the $7$-regular triangulation}
\author{Itai Benjamini}

\date{30.11.16}

\maketitle

\begin{abstract}
It is shown that self avoiding walk on the $7$-regular infinite planar triangulation has linear expected displacement.

\end{abstract}

\section{Introduction}

Consider $n$-steps {\em self avoiding walk}  (SAW) on the $7$-regular infinite planar triangulation.
The goal of this note it to prove
that the SAW is ballistic in expectation, i.e., that the expected distance from the starting point grows linearly with $n$.
This was already established, using counting arguments, for planar lattices of sufficiently high degree in \cite{MW}.
This argument fails for the $7$-regular triangulation.
In \cite{DH} it is proved that Euclidean SAW is sub-ballistic,
and see \cite{L} for ballisticity of SAW on Cayley  graphs with more than one end.
Recall that the $n$-steps SAW is the uniform measure on the set of all length $n$ self avoiding paths starting at a fixed root.
The next subsection contains a brief background on Gromov hyperbolic graphs,
for a thorough introduction to Gromov hyperbolic spaces see e.g. \cite{BI}.
For background on SAW see \cite{BDGS}.

\medskip

\begin{theorem}
Let $G$  be the $7$-regular infinite planar triangulation.
Then the expected displacement of $n$ steps SAW on $G$  is bigger than $cn$, where $c$ is a positive constant determined by $G$.
\end{theorem}

\noindent
{\sc Question:} Improve the theorem to an asymptotically almost sure ballisticity.

That is, there is $c >0$ so that the probability $n$-steps SAW reach distance $cn$ from the origin,
is going to $1$ with $n$.

No bounds on the number of self avoiding walks are used, unlike the argument in \cite{MW},
and therefore no bounds on the spectral gap are needed.
\medskip

The proof here works as is, to achieve expected ballisticity for every $1$-skeleton of any hyperbolic lattice in $\H^n$,
that is generated  by reflections of the fundamental domain. See \cite{V} for a survey on such lattices,
which exist  also in dimensions bigger than $2$, yet
there are no such of dimension bigger than $62$. We leave it as an open
{\em question} to adapt the proof to any finite volume lattice.

Another model in which the the same argument can give expected ballisticity,
is that of SAW defined directly on $\H^n$ as follows:
consider continuous $n$ steps random walk on $\H^n$, where the next step is chosen uniformly
and independently on the unit sphere around the current location.
Look at the restriction of this uniform product measure
to sequences in which the distance between any pair of vertices is bigger than $c$, for some
$0 < c < 1$, fixed.
\medskip

{\em Here is an outline of the proof}.  Given a self avoiding walk $\gamma$ in a space of constant negative curvature,
a constant portion of the vertices are at the boundary of the convex hull of the walk. After reflection the part of $\gamma$ from $\gamma(i)$ to $\gamma(n)$,
at a tangent of the convex hull at  $\gamma(i)$,
the inverse triangle inequality is satisfied at  $\gamma(i)$, that is $$ \dist(\gamma(1),\gamma(i)) + \dist(\gamma(i), \gamma(n)) ) \leq   \dist(\gamma'(1), \gamma'(n)) + \C, $$
where $\C < \infty$ depends only on the lattice, and $\gamma'$ denotes the reflected path.
We will see that the counting measure on the space of SAW of length $n$ is up to multiplicative constant, invariant with respect to reflections at a tangent of the convex hull, as these are order $n$ abundant.
Therefore a constant fraction of the vertices on a constant fraction  of the self avoiding walks, satisfy the  inverse triangle inequality.
Now a path in bounded degree Gromov  hyperbolic graph, with inverse triangle inequality holding for portion
of the vertices must be ballistic.
\medskip

Twenty years ago Oded Schramm briefly commented  over lunch,
that abundance of  reflections  will imply that hyperbolic SAW is ballistic.
Conceivably the argument here follows what Oded had in mind.

\subsection{Background}
For the probabilist reader, note that we use only well known basic properties of hyperbolic geometry.
We review them briefly here, see e.g.\cite{CFKP}.

Recall (see e.g.  \cite{BI}) that a graph is Gromov hyperbolic if for some $\delta < \infty$, all geodesic triangles are $\delta$-thin.
Where a geodesic triangle  is called $\delta$-thin if the
distance from any point on one of the geodesics to the union of the two other geodesics is bounded by $\delta$.

\begin{definition}
Let $G$ be an infinite connected graph, and let $C \geq 1$ and $K \geq 0$. A
path $\gamma$ of length $n$ between two vertices $x$ and $y$ is called a $(C, K)$-quasi-geodesic
finite path if for all $0 < i < j ≤ n$,
$$
j − i \leq Cd(\gamma(i), \gamma(j)) + K.
$$
\end{definition}

Similarly, we define $(C, K)$-quasi-geodesic infinite (or bi-infinite) paths. An
infinite (or a bi-infinite) path will simply be called a quasi-geodesic if it is
$(C, K)$-quasi-geodesic for some constants $C$ and $K$.

\begin{definition}
A bi-infinite path $\gamma$ in $G$ is called a Morse quasi-geodesic (resp.
Morse geodesic) if it is a quasi-geodesic (resp. a geodesic) and if it satisfies
the so-called Morse property: for all $C ≥ 1$ and $K > 0$, there exists $R$ such
that every $(C, K)$-quasi-geodesic joining two points of $\gamma$ remains inside the
$R$-neighborhood of $\gamma$.
\end{definition}

All geodesics in Gromov hyperbolic graphs are Morse (see e.g.  \cite{BI}).

The $7$-regular triangulation  as any co compact lattice in $\H^n$, is quasi-isometric to $\H^n$.
Consider the $7$-regular triangulation as a subset of the hyperbolic plane.
Since all two sided infinite geodedics are Morse, any quasi-geodesic is in a bounded neighborhood of a geodesic,
it follows that there exists a constant $C$ such that for every finite set of vertices $K$ in $G$, the image of the convex hull of $K$ in $G$ under the embedding into $\H^2$ has Hausdorff distance at most $C$ from the convex hull in $\H^2$ of the image of $K$ into $\H^2$ under the embedding. Where 

\begin{definition}
The convex hull of a set in a graph is the union of all the vertices on  geodesics between all pairs of vertices in the set.
\end{definition}

\begin{figure}[h]
\begin{minipage}[b]{1\linewidth}
\centering
\includegraphics[scale=0.5]{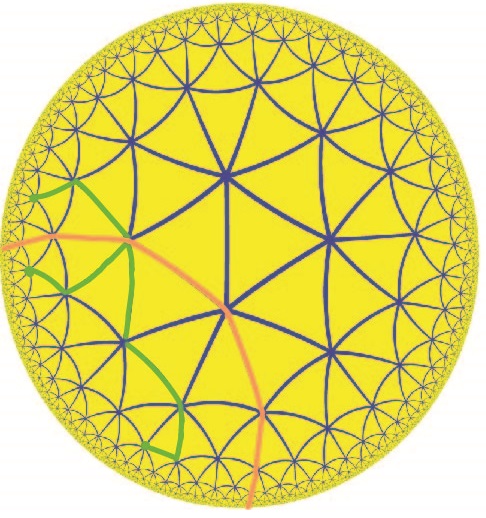}
\label{fig:7hyp}
\caption{A green path with  the last part reflected at the orange graph geodesic.}
\end{minipage}
\quad
\end{figure}

\section{Proof}

\begin{proof}

From now on $G$ denotes the $7$-regular infinite planar triangulation.
The following two lemmas are essentially from \cite{BE}. The convex hull is taken in the hyperbolic space.
As in Euclidean space, the convex hull is the intersection of all the half spaces containing the set.
As noted the convex hull in $G$ when imbedded in $\H^2$ has uniformly bounded hausdorff distance to the
convex hull in $\H^2$.

\begin{lemma}\label{convex}
Let $A$  a finite set of vertices in $G$. There exists $c>0$,
so that at least $c|A|$ vertices of $A$ are on  boundary of the convex hull of $A$.
\end{lemma}

More generally

\begin{lemma}\label{convex2}
Let $A$  a finite set of vertices in a hyperbolic lattice. Then there are $c>0, c'>0$ depending only on the lattice,
so that at least $c|A|$ vertices of $A$ are at distance at most $c'$ from the boundary of the convex hull of $A$.
\end{lemma}

{\sc Proof:} Let $\alpha$ be the minimal distance between distinct lattice points,
and let $A_\alpha = \{ x \in \H^n \, ; \, d(x, A) \leq \alpha \}$.
For any finite subset $A$ of the lattice, $c < Vol_n(A_\alpha) / |A| < C$.

According to \cite{BE}, the volume of the convex hull $K$ of $A_\alpha$  differs from $Vol_n(A_\alpha)$ by at most a constant factor.
Pick a random point $x \in K$. Since the lattice is non amenable, with constant probability, it lies in $A_\alpha$,
and its distance from $\partial K$ is at most a constant.
Therefore a constant fraction of $A_\alpha$ is of constant distance from $\partial K$, which implies
the same for $A$.    $\hfill \Diamond$

The following lemma follows directly from  the fundamental property of Gromov hyperbolic space,
that the distance outside a ball of radius $r$, between two geodesic rays starting distance $\delta$ apart,
diverges exponentially with $r$, see e.g. \cite{BI}.

\begin{lemma}\label{two convex sets }
Let $H$ be a $\delta$-hyperbolic graph.
Assume $A, B \subset H$,  are disjoint convex sets
with $d(A,B) := \inf_{x \in A, y \in B} d(x,y) \geq \delta$. Assume $x_0 \in A$ and $y_0 \in B$ realized the distance
$d(A,B)$, then the following holds:
There is $C$ depending only on $\delta$, such that for any $x \in A, y,   \in B$,
the geodesics between $x$ and $y$ intersects the ball
around $x_0$ of radius $C$.
\end{lemma}

\noindent
{\sc Proof of Lemma:}
In any Gromov hyperbolic space where all triangles are $\delta$-thin,
if the distance between  a pair of two sided infinite geodesics, is bigger  than $\delta$,
then the distance between the pair of geodesics outside a ball of radius $r$, grows at least exponentially in $r$.
\medskip

Separate the two convex sets by two such geodesics to finish.
(In higher dimension separate the two convex sets by geodesic hyperplanes, to get a similar result). $\hfill \Diamond$

Note that for the $7$-regular triangulation $\delta =1$ (exercise).

\begin{lemma} Let $A$ be a finite subset of vertices of $G$, and let $x_0 \in A$ belong to the boundary
of the graph convex hull, $K = conv(A)$. Then there exists a lattice-reflection $g$ such that $g(x_0) = x_0$ and $g(K) \cap K = \{ x_0 \}$.
\end{lemma}

\noindent
{\sc Proof of Lemma:} $G$ is obtained by consecutive reflections of a fixed geodesic triangle in $\H^2$.
In particular the orange line in the figure is a geodesic,
See also \cite{V}. Therefore if $X_0$ is on the convex hull, there is a graph geodesic via $x_0$, for which all $A$ is in
one half plane.       $\hfill \Diamond$

By the comments in the background section, if $x_0$ is within a bounded distance, {\em $d$} to the boundary of the convex hull, then there is a reflection for which the intersection is
contained in a ball of finite radius depending only on {\em $d$}.

\medskip

Denote by $\Lambda_n$ the set of $n$ steps SAW's.

$G$ is obtained by reflecting an heptagon (the dual graph of the triangulation) in $\H^2$, which is a fundamental domain for the seven regular lattice.
Pick the $7$  reflections along geodesics which extends the edges of the fundamental domain.
For each  vertex in the graph fix $7$ such reflections, and  order them arbitrarily.

Define  $R_i: \Lambda_n  \rightarrow \Lambda_n$ as follows.
Pick  the first  reflection, in the order on reflections at $\gamma(i)$,
which reflect the part of $\gamma$ from $\gamma(i)$ to $\gamma(n)$, along a tangent to the convex hull of the self avoiding walk at $\gamma(i)$, and apply it on $\gamma$.
If there is no such reflection at $\gamma(i)$, then $R_i$  doesn't change the walk.

For every $\gamma$, $|R_i^{-1} (\gamma)| \leq 7$.

\noindent
{\bf Remark:}
{\em (Regarding other lattices of  hyperbolic reflection groups)}
Note that for $1$-skeletons of other lattices arising from reflections,
when we only have  ~\ref{convex2}  and ~\ref{two convex sets }
with $c' < \infty $ and $\delta < \infty$,  a local modification of the  reflection at $i$, $R_i$ is needed:
reflect  the self avoiding walk along a geodesic tangent to a vertex $v$
of distance $\max\{ c', \delta\}$ away from $\gamma(i)$
and the convex hull, and then locally modify the reflected self avoiding walk
by connecting $\gamma(i)$ to $v$, by a $\max\{ c', \delta\}$ long path,
and remove $\max\{ c', \delta\}$ long segment from the end of the walk.


By lemma~\ref{convex}  uniformly over all walks in $\Lambda_n$, a constant portion of the vertices are
on the the boundary of the convex hull, so it is possible to apply a reflection at $\gamma(i)$ at a portion of the indexes,
and for any walk the total number of reflections is order $n$.

Let $\gamma$ be chosen according to the SAW counting measure and fix a constant $\C >0$. Denote by $\A_i$ the event that,
$$
\{\dist(\gamma(1),\gamma(i)) + \dist(\gamma(i), \gamma(n)) ) \leq   \dist(\gamma(1), \gamma(n)) + \C \}
$$

that is a  {\em $\C$-inverse triangle inequality} holds at $i$.

\begin{lemma}

There are $C < \infty$ and   $ c >0$ depending only on the graph $G$, so that for any $\gamma$ path of length $n$,  there are $cn$ indexes $i$, so that the event $\A_i$ occurs.

\end{lemma}

\noindent
{\sc Proof of Lemma:}
Fix $\gamma$ and pick $i$ uniformly between $0$ and $n$. With probability bigger than $c$, $\gamma(i)$ is at the boundary of the convex hull.
Observe that  for any $\gamma$ the number of reflections outside of the convex hull is proportional to the number of all the reflections, which is order $n$.
The probability of the event $\A_i$, is up to a constant factor
invariant with respect to $R_i$ and when $R_i$ reflects at $i$, the event above holds by lemma~\ref{two convex sets }. Therefore,
\begin{equation}
\begin{split}
 7  \times Prob(\gamma  \mbox{  satisfies  } \A_i) \geq Prob ( R_i \gamma  \mbox{  satisfies  } \A_i ) \\
 \geq Prob ( \gamma(i) \in \partial \mbox{convex hull of } \gamma ) \geq c.
\end{split}
\end{equation}
$\hfill \Diamond$

And thus,

\begin{equation}
\E |\{ i :  \C-\mbox{inverse triangle inequality holds at } i\}|  \geq cn.
\end{equation}

Since  geodesic triangles in Gromov hyperbolic graphs are $\delta$-thin, if $\C$-inverse triangle inequality holds at $i$, by Gromov hyperbolicity,
it has to be in a bounded neighborhood of the geodesic, depending only on the $\delta$-hyperbolicity and $\C$.

The number of vertices in $G$ in a bounded neighborhood of the geodesic interval
between $\gamma(1)$ and $\gamma(n)$ is bounded from above by $c' \dist(\gamma(1), \gamma(n))$, for some $c' < \infty$.
By (2.2) this neighborhood contains $cn$ vertices of $\gamma$ in expectation.

By self avoidance these $cn$ vertices are disjoint, thus the expected distance from $\gamma(1)$ to $\gamma(n)$ has to be
order $n$.

\end{proof}

\section{Comments on SAW on large graphs}

Recently Tom Hutchcroft  that SAW  on $\T_3 \times \Z$ is ballistic.
Where $\T_3$ is the three regular tree, \cite{H}.

SAW should be ballistic on any Cayley graph in which simple random walk is ballistic.
In particular all non amenable Cayley graphs.

Assume SAW is ballistic on $G \times \Z$, where $G$ is an infinite vertex transitive graph.
Show that SAW is ballistic on $G$.
Maybe \cite{DH} can be used to show sub ballisticity in the $\Z$ direction?

\medskip

We expect SAW to be ballistic on any infinite planar triangulation, in which all degrees are
at least $6$ (non positive curvature) and there is $r < \infty$,
so that every $r$-ball in the triangulation contains a vertex of degree at least $7$. Such triangulations are non amenable
(see \cite{ABH}) and for any finite set, a constant fraction of the vertices are at the boundary of the convex hull.

\medskip

Consider the uniform measure on connected subsets of size $n$ in a hyperbolic lattice,
or non amenable hyperbolic transitive graphs, containing a fix vertex $v_0$.

\noindent
{\sc Conjecture:} The scaling limit of uniform connected subsets is  Brownian tree.

See the survey \cite{LeG} for background on the Brownian tree.

Both when taking the ambient metric on the connected set, as an embedded subset, or the internal metric.
See \cite{CM} for a related result and definition of the scaling limit. We conjecture that  Brownian tree is also the scaling limit of
{\em self avoiding loop} in the ambient metric. In \cite{DGHM} it was shown that the probability Euclidean SAW makes a loop decays polynomially
to $0$. Show exponential decay for  the analogous event on the seven regular triangulations. This is a step towards ballsticity of the SAW.
Simple counting gives this for $d$-regular triangulation with large $d$.

\noindent
{\bf Acknowledgements:} I'm very grateful to Bo'az Klartag,
his comments clarified the argument and improved the presentation.
Thanks to Hugo  Duminil-Copin, for useful discussions and for
Tom Hutchcroft and Aran Raoufi for comments on a previous version.

\end{document}